\documentclass{amsart}
\usepackage{latexsym,amsxtra,amscd,ifthen}
\usepackage{amsfonts}
\usepackage{verbatim}
\usepackage{amsmath}
\usepackage{amsthm}
\usepackage{amssymb}

\numberwithin{equation}{section}

\theoremstyle{plain}
\newtheorem{theorem}{Theorem}[section]
\newtheorem{lemma}[theorem]{Lemma}
\newtheorem{proposition}[theorem]{Proposition}
\newtheorem{corollary}[theorem]{Corollary}
\newtheorem{conjecture}[theorem]{Conjecture}
\newtheorem{construction}[theorem]{Construction}
\newtheorem{question}[theorem]{Question}
\newtheorem{problem}[theorem]{Problem}

\theoremstyle{definition}
\newtheorem{definition}[theorem]{Definition}

\newtheorem*{remark*}{Remark}

\newtheorem{subsec}[theorem]{}

\newcommand{\ZZ}{{\mathbb Z}}

\newcommand\calO{{\mathcal O}}

\newcommand\eps{\varepsilon}
\newcommand\gfrak{{\mathfrak g}}

\newcommand\Ug{U(\gfrak)}
\newcommand\Uqg{U_q(\gfrak)}
\newcommand\OG{\calO(G)}
\newcommand{\Oq}{\calO_q}
\newcommand{\OqG}{\Oq(G)}

\DeclareMathOperator{\gldim}{gl{.}dim}
\DeclareMathOperator{\rgldim}{r{.}gl{.}dim}
\DeclareMathOperator{\lgldim}{l{.}gl{.}dim}
\DeclareMathOperator{\pdim}{p{.}dim}
\DeclareMathOperator{\Ext}{Ext}
\DeclareMathOperator{\GKdim}{GKdim}

\newcommand{\id}{\operatorname{id}}

\newcommand\chr{\operatorname{char}}
\newcommand\bull{$\bullet$\ }
\newcommand\Znonneg{\ZZ_{\ge0}}

\newcommand\kx{k^{\times}}
\DeclareMathOperator{\pideg}{PIdeg}

\begin{document}

\title[Noetherian Hopf algebras]
{Noetherian Hopf algebras}

\author{K. R. Goodearl}

\address{Department of Mathematics, 
University of California at Santa Barbara,
Santa Barbara, CA 93106, USA}

\email{goodearl@math.ucsb.edu} 

\begin{abstract}
A brief survey of some aspects of noetherian Hopf algebras is given, concentrating on structure, homology, and classification, and accompanied by a panoply of open problems.
\end{abstract}




\dedicatory{Dedicated to Kenny Brown and Toby Stafford on the occasion of their sixtieth birthdays.}

\maketitle

\setcounter{section}{-1}
\section{Introduction}

At a 1997 conference in Seattle, Kenny Brown promoted ``noetherian Hopf algebras" as a subject to be studied in its own right, surveyed the status of this subject, and proposed a number of problems and conjectures \cite{Br2}. A decade later, he updated the status in the survey \cite{Br3}. Our purpose here is to give a brief further update, and especially to highlight a good selection of open problems. Many of these are shamelessly borrowed from \cite{Br2}, \cite{Br3}. The background is well laid out in those surveys, which readers are encouraged to look at in conjunction with the present paper.

Motivation for studying noetherian Hopf algebras comes in part from the fact that the examples known for many years share a number of characteristics. Further, it is important to understand those that appeared with the rise of quantum groups. From the most general perspective, as put forward by Drinfel'd \cite[p.~800]{Dri}, the notions of ``Hopf algebra" and ``quantum group" are equivalent. Specializing to affine noetherian Hopf algebras can thus be viewed as developing the theory of ``quantum affine algebraic groups". 

One might ask why it is reasonable to attack noetherian Hopf algebras as a class, i.e., why is this not a hopelessly large and amorphous mixed bag of objects, like the class of all rings? One response is that the existence of a Hopf algebra structure on an algebra imposes a lot of rigidity -- for instance, the representations form a monoidal category. Another is to point to the extensive theory which has been developed for finite dimensional Hopf algebras. The existence of this theory also prompts one to seek suitably modified analogs in the infinite dimensional, noetherian case. On the other hand, the extent of the finite dimensional theory rules against attempting to subsume this in a noetherian theory. Taking existing examples together with the principles of quantum groups into account, and recalling that the coordinate rings of connected algebraic groups are integral domains, one is prompted to focus the study of noetherian Hopf algebras primarily on those which are domains, or at least prime rings. This has the advantage of avoiding all finite dimensional Hopf algebras except for the trivial $1$-dimensional one (although many noetherian Hopf algebras have nontrivial finite dimensional Hopf algebra quotients).

Leaving finite dimensional Hopf algebras aside, the starting point for noetherian Hopf algebra theory consists of five basic, long-known classes of examples. The Hopf structures in each case are so well known that we do not describe them here, but just refer to \cite{Br2}, \cite{Br3}, for instance.
\begin{enumerate}
\item[\bull] The enveloping algebra $\Ug$ of any finite dimensional Lie algebra $\gfrak$. (While all enveloping algebras are Hopf algebras, it is not known whether there exists any infinite dimensional Lie algebra whose enveloping algebra is noetherian.) 
\item[\bull] The group algebra $k\Gamma$ of any polycyclic-by-finite group $\Gamma$, over a field $k$. (While all group algebras are Hopf algebras, it is not known whether there exists any non-[polycyclic-by-finite] group whose group algebra is noetherian.)
\item[\bull] The coordinate ring $\OG$ of any affine algebraic group $G$.
\item[\bull] The quantized enveloping algebra $\Uqg$ of any finite dimensional semisimple complex Lie algebra $\gfrak$, where $q$ is a nonzero scalar which may or may not be a root of unity. (Descriptions of these Hopf algebras are given, for instance, in \cite[Chapter I.6]{BG}. One can find variants, multiparameter versions, and cocycle twists of these in the literature.)
\item[\bull] The quantized coordinate ring $\OqG$ of any connected semisimple complex algebraic group $G$, where $q$ again is an arbitrary nonzero scalar. (Descriptions are given, for instance, in \cite[Chapter I.7]{BG}. Variants, multiparameter versions, and cocycle twists occur in the literature.)
\end{enumerate}
An immediate guiding question, of course, is to ask how representative of noetherian Hopf algebras in general these examples might be.

The following discussion concentrates on three aspects of the theory: structure, homology, and classification. The main lines of what is known will be sketched, accompanied by many questions.

\begin{subsec}
Fix a base field $k$. Throughout, all vector spaces, algebras, tensor products, etc.~are taken over $k$. For any Hopf algebra, we denote the multiplication, unit map, comultiplication, counit, and antipode by $m$, $u$, $\Delta$, $\eps$, and $S$, respectively. A recommended reference for general Hopf algebra theory is \cite{Mon}; an outline of the basic concepts can be found in \cite[Appendix I.9]{BG}.
\end{subsec}

\section{Structure}
\label{yysec1}

It is obviously too much to ask for a specific structure theory covering arbitrary noetherian Hopf algebras, but a number of general issues raise their heads. Since all the important noetherian examples are also affine (i.e., finitely generated as algebras), we first ask about the relationship between these conditions.

\section*{Noetherian versus affine}

\begin{question}
How are the noetherian and affine conditions related for Hopf algebras?
\end{question}

To support the reasonableness of this question, we offer the following results of Molnar \cite[Proposition and following remarks]{Mol} and Liu-Zhang \cite[Corollary]{LZ}.

\begin{theorem}
{\rm [Molnar] (a)} A commutative Hopf algebra is noetherian if and only if it is affine.

{\rm(b)} Any cocommutative noetherian Hopf algebra is affine.
\end{theorem}

\begin{theorem}
{\rm [Liu-Zhang]} A Hopf algebra is artinian if and only if it is finite dimensional {\rm(}in which case it is also affine\/{\rm)}.
\end{theorem}

In general, affine Hopf algebras need not be noetherian, as shown by the group algebra $k\Gamma$ of a nonabelian free group $\Gamma$ of finite rank. Whether the converse might hold is an open question, raised by Wu and Zhang \cite[Question 5.1]{WZ2}. In view of Molnar's theorem, one might still ask for the implication ``affine$\implies$noetherian" in situations ``close to commutative". Brown has raised this question for Hopf algebras satisfying a polynomial identity \cite[Question C(iii)]{Br3}. 

\begin{question}
{\rm (a) [Wu-Zhang]} Is every noetherian Hopf algebra affine?

{\rm (b) [Brown]} Is every affine PI Hopf algebra noetherian?
\end{question}

Among PI algebras, affineness often comes in the form of finiteness over the center. The latter property does not always hold, however, even for affine, noetherian, PI Hopf algebras, as shown by a construction of Gelaki and Letzter \cite[Remark 3.9]{GeLe} (see \cite[p.~10, footnote]{Br3}). Brown specialized the question to the semiprime case \cite[Question C(i)]{Br3}:

\begin{question}
{\rm [Brown]} Is every semiprime noetherian PI Hopf algebra module-finite over its center? What if it is also assumed to be affine?
\end{question}

\section*{The antipode}

Recall that in the definition of a Hopf algebra, the antipode is only required to be invertible with respect to the convolution product, not necessarily as a map. All that follows from the axioms is that $S$ must be an algebra anti-endomorphism. Bijectivity of the antipode is important in many analyses, and appears in all known noetherian examples, so this condition deserves some focus. To start with the negative side: Takeuchi constructed Hopf algebras with non-bijective antipodes in \cite[Theorem 11]{Tak1}, \cite{Tak2}. On the other hand, the situation is fine in the finite dimensional case, due to the following theorem of Radford \cite[Theorem 1]{Rad}:

\begin{theorem}
{\rm [Radford]} The antipode of any finite dimensional Hopf algebra is bijective; in fact, it has finite order {\rm(}i.e., some power of $S$ is the identity\/{\rm)}.
\end{theorem}

The strongest conclusion in this direction is obtained for commutative or cocommutative Hopf algebras \cite[Proposition 4.0.1(6)]{Swe}, \cite[Corollary 1.5.12]{Mon}.

\begin{proposition}
In any commutative or cocommutative Hopf algebra, $S^2= \id$.
\end{proposition}

Skryabin has established bijectivity of the antipode for many noetherian Hopf algebras \cite[Corollaries 1, 2]{Skr}, and made a general conjecture \cite[Conjecture]{Skr}.

\begin{theorem}
{\rm [Skryabin] (a)} The antipode of any noetherian Hopf algebra is injective.

{\rm (b)} Let $H$ be a Hopf algebra which is either semiprime noetherian or affine PI. Then $S$ is bijective.
\end{theorem}

\begin{conjecture}
{\rm [Skryabin]} The antipode of every noetherian Hopf algebra is bijective.
\end{conjecture}

Natural related questions, taking Radford's theorem into account, were raised by Brown and Zhang \cite[Questions H, I]{Br3}, \cite[Question 6.2]{BZ1}.

\begin{question}
{\rm [Brown, Brown-Zhang]} If $H$ is an affine noetherian PI Hopf algebra, does $S$ have finite order? Is some {\rm(}even{\rm)} power of $S$ an inner automorphism of $H$?
\end{question}

In a non-PI noetherian Hopf algebra, the antipode need not have an inner power. For example, let $H= \Oq(SL_2(k))$, which is an affine noetherian Hopf algebra. If $q\in\kx$ is not a root of unity, then all powers of $S$ are outer.

Some positive cases for the above question were obtained by Brown and Zhang \cite[Proposition 6.2 and following comments]{BZ1}.

\begin{theorem}
{\rm [Brown-Zhang]} Let $H$ be an affine noetherian PI Hopf algebra. If $H$ either has finite global dimension or is module-finite over its center, then some power of $S$ is inner.
\end{theorem}

\section*{Grouplike and skew primitive elements}

Two of the most important models for noetherian Hopf algebras are group algebras and enveloping algebras. These are distinguished by the behavior of the comultiplication: in $k\Gamma$, we have $\Delta(g)= g\otimes g$ for all $g\in \Gamma$, whereas in $\Ug$, we have $\Delta(x)= x\otimes1+ 1\otimes x$ for all $x\in \gfrak$. Conversely, one can collect elements of these types in any Hopf algebra $H$ to obtain a group or a Lie algebra. More precisely, the set of all grouplike elements of $H$ (see below) is a subgroup of the group of units of $H$, while the set of all primitive elements is a Lie algebra with respect to the additive commutator. We recall the definitions.

\begin{definition}
Let $H$ be a Hopf algebra.

An element $g\in H$ is \emph{grouplike} provided $g\ne0$ and $\Delta(g)= g\otimes g$. (It then follows from the counit axiom that $\eps(g)=1$.)

An element $x\in H$ is \emph{primitive} provided $\Delta(x)= x\otimes1+ 1\otimes x$. (The counit axiom then implies $\eps(x)=0$.) More generally, $x$ is \emph{$(g,h)$-skew primitive}, where $g$ and $h$ are grouplike elements of $H$, provided $\Delta(x)= x\otimes g+ h\otimes x$. (Again, it follows that $\eps(x)=0$.)
\end{definition}

Although it might seem overly myopic to concentrate on such elements, they play a central enough role to be worthy of focus. The following easy lemma gives the basic relationships with the antipode and a useful consequence.

\begin{lemma}  \label{Sform}
Let $H$ be a Hopf algebra.

{\rm (a)} If $g\in H$ is grouplike, then $g$ is invertible, $S(g)= g^{-1}$, and $g^{-1}$ is grouplike.

{\rm (b)} If $x\in H$ is $(g,h)$-skew primitive, for some grouplike $g,h\in H$, then $S(x)= -h^{-1}xg^{-1}$.
\end{lemma}

It follows from Lemma \ref{Sform} that if a Hopf algebra $H$ is generated {\rm(}as an algebra\/{\rm)} by its grouplike and skew primitive elements, then $S$ is surjective. In fact, $S$ must be bijective in this situation (Corollary \ref{corSbijec}).

\begin{question}
Which noetherian Hopf algebras are generated by their grouplike and skew primitive elements?
\end{question}

For a negative example, consider $\calO(SL_2(k))$. This noetherian Hopf algebra has no nonzero skew primitive elements, and no grouplike elements except $1$.

Clues to the question of generation by grouplikes and skew primitives can be obtained from the coalgebra structure of a Hopf algebra, as follows.

\begin{definition}
A \emph{subcoalgebra} of a Hopf algebra $H$ is any linear subspace $C\subseteq H$ such that $\Delta(C) \subseteq C\otimes C$. It is a \emph{simple subcoalgebra} provided $C\ne \{0\}$ and the only subcoalgebras contained in $C$ are $\{0\}$ and $C$. The Hopf algebra $H$ is called \emph{pointed} if every simple subcoalgebra of $H$ is $1$-dimensional.
\end{definition}

The connection with the above discussion is that the $1$-dimensional subcoalgebras of a Hopf algebra $H$ are precisely the subspaces $kg$ for grouplike $g\in H$.

For instance, over an algebraically closed base field, cocommutative Hopf algebras are pointed \cite[Lemma 8.0.1(c)]{Swe}. More specific examples of pointed Hopf algebras are $k\Gamma$, $\Ug$, and $\Uqg$. On the other hand, $\calO(SL_2(k))$ is not pointed, since it has a $4$-dimensional simple subcoalgebra, namely the subspace spanned by the matrix entry functions $X_{11}$, $X_{12}$, $X_{21}$, $X_{22}$. Since examples of this type only seem to appear in dimension $3$ and higher, Brown and Zhang raised the following question \cite[\S0.5]{BZ2}.

\begin{question}  \label{qnGK2pointed}
{\rm [Brown-Zhang]} Let $H$ be a prime affine noetherian Hopf algebra. If $\GKdim(H)\le 2$, is $H$ pointed?
\end{question}

The most useful condition implying pointedness is the existence of enough group\-like and skew primitive elements, as the following lemma shows. It is a corollary of \cite[Lemma 5.5.1]{Mon}.

\begin{lemma}
If a Hopf algebra $H$ is generated {\rm(}as an algebra\/{\rm)} by its grouplike and skew primitive elements, then $H$ is pointed.
\end{lemma}

The converse fails, even in the noetherian case. Here are two examples, taken from \cite[Examples 5.11, 5.12]{AnS2000}; a third will appear following Question \ref{WZZqn}. For the first, assume $\chr(k)=p>0$ and view the polynomial ring $k[x]$ as the enveloping algebra of the $1$-dimensional Lie algebra $kx$. Then $k[x]$ is a Hopf algebra, with $x$ primitive. Due to characteristic $p$, the ideal $\langle x^{p^2}\rangle$ is a Hopf ideal of $k[x]$, and so $k[x]/\langle x^{p^2}\rangle$ is a finite dimensional Hopf algebra. Its dual, $\bigl( k[x]/\langle x^{p^2}\rangle \bigr){}^*$, is a finite dimensional (and thus noetherian) pointed Hopf algebra which is not generated by its grouplikes and skew primitives.

For the second example, take $\chr(k)=0$, let $q\in\kx$ be a primitive $\ell$-th root of unity for some $\ell>1$ (e.g., take $q=-1$), and let $H$ be a vector space with basis $\{x_n\mid n\in\Znonneg \}$. Then $H$ can be made into a Hopf algebra in which 
$$x_mx_n= \begin{pmatrix}m+n\\  n\end{pmatrix}_q x_{m+n} \qquad\text{and}\qquad \Delta(x_n)= \sum_{i=0}^n x_i \otimes x_{n-i}$$
for all $n$. It is noetherian and pointed, but not generated by its grouplikes and skew primitives.

That a characteristic zero example must be infinite dimensional is the content of a conjecture of Andruskiewitsch and Schneider \cite[Conjecture 1.4]{AnS2000}.

\begin{conjecture}
{\rm [Andruskiewitsch-Schneider]} Assume $k$ is algebraically closed of characteristic zero. Then any finite dimensional pointed Hopf algebra over $k$ is generated by its grouplikes and skew primitives.
\end{conjecture}

One case of this conjecture has been established by Angiono \cite[Theorem 2]{Ang}:

\begin{theorem}
{\rm [Angiono]} Let $H$ be a finite dimensional pointed Hopf algebra over an algebraically closed field of characteristic zero. If the group of grouplike elements of $H$ is abelian, then $H$ is generated by its grouplikes and skew primitives.
\end{theorem}

A sharpening of Question \ref{qnGK2pointed}, in the case of a domain, is suggested by comments of Wang, Zhang and Zhuang \cite[Comments following Corollary 0.2]{WZZ1}.

\begin{question}  \label{WZZqn}
{\rm [Wang-Zhang-Zhuang]} Let $H$ be an affine noetherian Hopf algebra domain with $\GKdim(H) \le 2$. Is $H$ generated by its grouplikes and skew primitives?
\end{question}

The restriction to $\GKdim(H)\le 2$ in Question \ref{WZZqn} is necessitated by the example of $\calO(SL_2(k))$. In fact, even in the pointed case, an affine noetherian Hopf algebra domain with Gelfand-Kirillov dimension greater than $2$ need not be generated by its grouplikes and skew primitives. The following example was shown to us by J.J. Zhang: take $H=\calO(G)$ where $G$ is the group of unipotent upper triangular $3\times3$ matrices over an infinite field $k$.  Then $H$ is a polynomial ring $k[x,y,z]$, with $x$ and $y$ primitive while $\Delta(z)= z\otimes1+ x\otimes y+ 1\otimes z$. The only grouplike element of $H$ is $1$, and the only (skew) primitive elements are the linear combinations of $x$ and $y$. Thus, $H$ is an affine noetherian Hopf algebra domain not generated by its grouplikes and skew primitives. An application of \cite[Lemma 5.5.1]{Mon} shows that all simple subcoalgebras of $H$ are contained in $k[x]$, and therefore $H$ is pointed. 

To end the section, we return to bijectivity of antipodes. The following ``folklore" result is given in \cite[Corollary 5.2.11]{Mon}.

\begin{proposition}
Let $H$ be a Hopf algebra. If all simple subcoalgebras of $H$ are cocommutative, then $S$ is bijective.

In particular, the antipode of any pointed Hopf algebra is bijective.
\end{proposition}

\begin{corollary}  \label{corSbijec}
If a Hopf algebra $H$ is generated by its grouplike and skew primitive elements, then its antipode is bijective.
\end{corollary}

\section{Homological conditions}
\label{yysec2}

All known noetherian Hopf algebras enjoy strong homological properties; we discuss these next. Recall that for noncommutative noetherian rings, finite global dimension alone is not a strong enough property to be very useful. Good upgrades include the versions of regularity introduced by Auslander and Artin-Schelter. Similarly, the Auslander and Artin-Schelter versions of the Gorenstein condition are the most useful in place of finite injective dimension. Good companions for these properties are the Cohen-Macaulay conditions, with respect to Krull or Gelfand-Kirillov dimension. Definitions for all these conditions are given in myriad sources; e.g., \cite[Appendix I.15]{BG}.

A key motivating result is the following theorem of Larson and Sweedler in the finite dimensional setting \cite[Remark, p.~85]{LS}, \cite[Theorem 2.1.3]{Mon}.

\begin{theorem}
{\rm [Larson-Sweedler]} Any finite dimensional Hopf algebra is a Frobenius algebra.
\end{theorem}

\begin{corollary}  \label{fdssiff}
Any finite dimensional Hopf algebra $H$ is self-injective. Consequently, $\gldim H<\infty$ if and only if $H$ is semisimple.
\end{corollary}

The right (or left) global dimension of a Hopf algebra $H$ is determined by the projective dimension of a single $H$-module, as follows. This result was known in several settings with additional hypotheses; the general result was observed by Lorenz and Lorenz \cite[Section 2.4]{LL}. Here $k_H$ and $_Hk$ denote the trivial $1$-dimensional right and left $H$-modules, respectively.

\begin{theorem}
$\rgldim H = \pdim k_H$ and $\lgldim H = \pdim {}_Hk$ for any Hopf algebra $H$.
\end{theorem}

Among our standard examples of noetherian Hopf algebras, $\Ug$, $\OG$, $\Uqg$, and $\OqG$ are Auslander-regular and Cohen-Macaulay \cite[Theorem B]{BG1}, \cite[Section 6]{BZ1}. If $\Gamma$ is a polycyclic-by-finite group, then $k\Gamma$ is Auslander-Gorenstein, but not necessarily Cohen-Macaulay \cite[Theorem 6.7, Remark 6.7(b)]{BZ1}.

For many purposes, it suffices to consider the Auslander versions of the above conditions, since they imply the Artin-Schelter versions by the following theorem of Brown and Zhang \cite[Lemma 6.1]{BZ1}.

\begin{theorem}  \label{AustoAS}
{\rm [Brown-Zhang]} Let $H$ be a noetherian Hopf algebra. If $H$ is Cohen-Macaulay and Auslander-regular {\rm(}respectively, Auslander-Gorenstein\/{\rm)}, then it is also Artin-Schelter-regular {\rm(}respectively, Artin-Schelter-Gorenstein\/{\rm)}.
\end{theorem}

Many authors have raised the question of extending the Larson-Sweedler theorem -- more precisely, the first part of Corollary \ref{fdssiff} -- to infinite dimensional noetherian Hopf algebras \cite[\S1.15]{BG1}, \cite[Question A]{Br2}, \cite[Question 0.3]{WZ1}, \cite[Question 5.2, Remark 5.9]{WZ2}, \cite[Question E]{Br3}.

\begin{question}
{\rm [Brown-Goodearl, Wu-Zhang]} Does every noetherian Hopf algebra $H$ have finite injective dimension? Is $H$ Auslander-Gorenstein? What if $H$ is also affine?
\end{question}

Wu-Zhang and Brown also posed the question whether a noetherian Hopf algebra in characteristic zero which is semiprime, or even a domain, must have finite global dimension \cite[Question 0.4]{WZ1}, \cite[Question K]{Br3}. The answer is negative, as shown by Goodearl and Zhang \cite[Remark 1.7]{GZ2}.

Several positive results in the homological direction have been proved, culminating in the following theorems of Brown-Goodearl \cite[Corollary 1.8]{BG1}, Wu-Zhang \cite[Theorem 0.1]{WZ1}, \cite[Theorems 0.1, 0.2]{WZ2}, and Lu-Wu-Zhang \cite[Theorem 0.4]{LWZ2}.

\begin{theorem}  \label{PIAGCM}
Let $H$ be a Hopf algebra.

{\rm(a)} {\rm [Brown-Goodearl, Wu-Zhang]} If $H$ is affine, noetherian, and PI, then it is Auslander-Gorenstein, Cohen-Macaulay, and semiprime.

{\rm(b)} {\rm [Wu-Zhang]} If $H$ is module-finite over its center, $Z(H)$ is affine, $S^2=\id_H$, and $\chr k=0$, then $H$ is Auslander-regular, Cohen-Macaulay, and semiprime.
\end{theorem}

\begin{theorem}
{\rm [Lu-Wu-Zhang]} Let $H$ be a nonnegatively filtered noetherian Hopf algebra whose associated graded algebra is connected graded with enough normal elements. Then $H$ is Auslander-Gorenstein and Cohen-Macaulay, and its antipode is bijective. If it has finite global dimension, then it is semiprime.
\end{theorem}

A nonnegatively graded algebra $A= \bigoplus_{n\ge0} A_n$ is \emph{connected graded} provided $A_0=k$, and it has \emph{enough normal elements} provided that for each homogeneous prime ideal $P \ne \bigoplus_{n>0} A_n$, the quotient $A/P$ contains a homogeneous normal element of positive degree.

We close this section with the following questions of Brown and Goodearl \cite[\S1.9]{BG1}:

\begin{question}
{\rm [Brown-Goodearl]} Let $H$ be a noetherian Hopf algebra with finite global dimension. Is $H$ Auslander-regular? Is it semiprime?
\end{question}

\section*{Integrals}

Regularity conditions in Hopf algebras are closely tied to concepts of integrals. There are left and right hand versions, with symmetric definitions; we concentrate on the latter.

\begin{definition}  \label{classicalint}
The set of \emph{right integrals} in a Hopf algebra $H$ is the set
$$\int^r_H := \{ x\in H \mid xa=\eps(a)x \text{\;for all\;} a\in H \}.$$
\end{definition}

Larson and Sweedler characterized semisimplicity (equivalently, finite global dimension) in finite dimensional Hopf algebras by means of integrals \cite[Proposition 3]{LS}. 

\begin{theorem}  \label{LSfdss}
{\rm [Larslon-Sweedler]} Let $H$ be a finite dimensional Hopf algebra. Then $H$ is semisimple if and only if $\eps\left( \int^r_H \right) \ne \{0\}$.
\end{theorem} 

Since infinite dimensional Hopf algebras often have no nonzero integrals, Lu, Wu, and Zhang introduced the following homological version \cite[Definition 1.1]{LWZ1}. They used it to characterize finite global dimension in certain cases \cite[Theorem 0.1]{LWZ1}.

\begin{definition}
Let $H$ be an Artin-Schelter-Gorenstein Hopf algebra with injective dimension $d$. In particular, this requires that $\Ext^i_H(k_H,H_H)$ vanishes for $i\ne d$, while $\Ext^d_H(k_H,H_H) $ is $1$-dimensional.
The \emph{right homological integral} of $H$ is
$$\int^r_H := \Ext^d_H(k_H,H_H),$$
which has a natural $H$-$H$-bimodule structure.
(In case $H$ is finite dimensional, there is a natural isomorphism between this version of $\int^r_H$ and the classical one defined as in Definition \ref{classicalint} \cite[p.~4948]{LWZ1}.) 
\end{definition}

Recall from Theorems \ref{PIAGCM}(a) and \ref{AustoAS} that affine noetherian PI Hopf algebras are Artin-Schelter-Gorenstein.

\begin{theorem}  \label{LWZfingldim}
{\rm [Lu-Wu-Zhang]} Let $H$ be an affine noetherian PI Hopf algebra with injective dimension $d$. Then $\gldim H <\infty$ if and only if
\begin{enumerate}
\item $\eps^*: \Ext^d_H\bigl(\int^r_H,{}_HH\bigr) \longrightarrow \Ext^d_H\bigl( \int^r_H,{}_Hk\bigr)$ is an isomorphism.
\item $\Ext^d_H(T,{}_Hk) =0$ for all simple left $H$-modules $T\not\cong \int^r_H$.
\end{enumerate}
\end{theorem}

A natural question, buttressed by Theorem \ref{LSfdss}, is whether condition (b) is needed in this theorem \cite[Question 3.6]{LWZ1}.

\begin{question}
{\rm [Lu-Wu-Zhang]} In the context of Theorem {\rm\ref{LWZfingldim}}, does condition {\rm(a)} alone imply $\gldim H <\infty$?
\end{question}

One further type of integral is the following, distinguished from those in Definition \ref{classicalint} by a small change in terminology. 

\begin{definition} A \emph{right integral \underbar{on}} a Hopf algebra $H$ is any functional $t\in H^*$ such that $(t\otimes f)\circ\Delta= f(1)t$ for all $f\in H^*$. (If $H$ is finite dimensional, this is the same as $t$ being a right integral in the dual Hopf algebra $H^*$, that is, $t\in \int^r_{H^*}$ in the sense of Definition \ref{classicalint}.)
\end{definition}

Integrals on Hopf algebras bring us back to antipodes once again, via the following result of Sweedler \cite[Corollary 5.1.7]{Swe} and Radford \cite[Proposition 2]{Rad2}:

\begin{theorem}
{\rm [Sweedler, Radford]} Let $H$ be a Hopf algebra. If there exists a nonzero right {\rm(}or left\/{\rm)} integral on $H$, then the antipode of $H$ is bijective.
\end{theorem}

\section{Classification}
\label{yysec3}

The rigidity of the structure of a Hopf algebra suggests that Hopf algebras are not thick on the ground, especially if there is some kind of limitation on their size. There is a big ongoing project to classify finite dimensional Hopf algebras, and there are plenty of these. Certainly, there is no hope of subsuming that project in any classification of noetherian Hopf algebras. From the noetherian viewpoint, however, vector space dimension is not a very important measure of size. In that arena, Krull dimension and Gelfand-Kirillov dimension take precedence. Furthermore, from the noncommutative viewpoint, the major algebras of interest are prime or even domains. Thus, it is reasonable to see how far prime noetherian Hopf algebras of low Gelfand-Kirillov dimension can be classified. The case of Gelfand-Kirillov dimension zero can be dismissed, since the only prime finite dimensional Hopf algebra over $k$ is $k$ itself.

\section*{Gelfand-Kirillov dimension $1$}

Lu, Wu, and Zhang initiated the program of classifying affine noetherian Hopf algebras of Gelfand-Kirillov dimension one in \cite[Section 7]{LWZ1}. This was carried forward by Brown and Zhang in \cite{BZ2}  under the following hypotheses:
\begin{enumerate}
\item[(H1)] $H$ is a prime, affine, noetherian Hopf algebra, the base field $k$ is algebraically closed of characteristic zero, $\gldim H<\infty$, and $\GKdim H =1$.
\end{enumerate}
By the Small-Stafford-Warfield Theorem \cite[Theorem]{SSW}, $H$ is then module-finite over its center; in particular, $H$ satisfies a polynomial identity. Further, it follows from Theorem \ref{PIAGCM}(a) that, in fact, $\gldim H = 1$.

There are three immediate examples satisfying (H1): 
\begin{enumerate}
\item[\bull] $\Ug$, where $\dim\gfrak=1$.
\item[\bull] $k\Gamma$, where $\Gamma=\ZZ$ or $\Gamma= \langle x,g\mid g^2=1,\, gxg^{-1}= x^{-1}\rangle$.
\end{enumerate}
The first two of these may also be presented as the coordinate rings of the algebraic groups $(k,{+})$ and $(\kx,{\cdot})$. There are no other coordinate rings to consider, since the two groups mentioned are the only connected one-dimensional algebraic groups over $k$. 

Two additional families are known to appear under (H1). One was introduced in \cite[Example 2.7]{LWZ1}; the other \cite[Section 3.4]{BZ2} generalizes a family constructed by Liu \cite[Section 2]{Liu}.
\begin{enumerate}
\item[\bull] $H(n,t,\xi)$, generated by a grouplike element $g$ and a $(g^t,1)$-skew primitive element $x$ satisfying $g^n=1$ and $xg=\xi gx$, where $n\ge2$ and $0\le t\le n-1$ are integers and $\xi\in\kx$ is a primitive $n$\,th root of unity.
\item[\bull] $B(n,w,\xi)$, generated by commuting grouplike elements $x^{\pm1}$, $g^{\pm1}$ and a $(g,1)$-skew primitive element $y$ satisfying $xy=yx$, $yg= \xi gy$, and $y^n= 1-x^w= 1-g^n$, where $n\ge2$ and $w\ge1$ are integers and $\xi\in\kx$ is a primitive $n$\,th root of unity.
\end{enumerate}

Brown and Zhang proved that in the case of prime PI-degree (and somewhat more generally), these examples constitute a classification of (H1) \cite[Theorem 0.5]{BZ2}. They also raised natural companion questions \cite[Questions 7.1, 7.2, 7.3C]{BZ2}.

\begin{theorem}
{\rm [Brown-Zhang]} If $H$ satisfies {\rm(H1)} and $\pideg H$ is prime, then $H$ is isomorphic to one of the examples described above.
\end{theorem}

\begin{question}
{\rm [Brown-Zhang] (a)} Do the above examples exhaust all Hopf algebras satisfying {\rm(H1)} with non-prime PI-degrees?

{\rm (b)} Is there a similar classification if the condition ``$\gldim H<\infty$" is removed from {\rm(H1)} and/or ``prime" is weakened to ``semiprime"?
\end{question}

\section*{Gelfand-Kirillov dimension $2$}

In \cite{GZ2}, Goodearl and Zhang moved to Gelfand-Kirillov dimension $2$, with slightly different hypotheses than in (H1), namely, restricting from prime rings to domains but allowing infinite global dimension. As it turns out, the classification results at this level for affine Hopf algebras are exactly the same as for noetherian Hopf algebras, so both can be stated together. The specific hypotheses are
\begin{enumerate}
\item[(H2)] $H$ is a Hopf algebra domain which is either affine or noetherian, the base field $k$ is algebraically closed of characteristic zero, and $\GKdim H=2$.
\end{enumerate}

There are several immediate examples satisfying (H2):
\begin{enumerate}
\item[\bull] $\Ug$, where $\dim\gfrak = 2$.
\item[\bull] $k\Gamma$, where $\Gamma= \ZZ^2$ or $\Gamma= \ZZ\rtimes\ZZ$ and in the semidirect product, $\ZZ$ acts on itself by the rule $m.n= (-1)^mn$.
\item[\bull] $\OG$, where $G= (k,{+})\rtimes(\kx,{\cdot})$ and $(\kx,{\cdot})$ acts on $(k,{+})$ by the rule $b.a= b^na$, for some $n\in\Znonneg$.
\end{enumerate}
The two group algebras mentioned may also be presented as the coordinate rings of $(k,{+})^2$ and $(\kx,{\cdot})^2$. These two groups, together with the ones mentioned in the last item, are the only connected two-dimensional algebraic groups over $k$.

Three additional families of examples satisfying (H2) were constructed in \cite[Section 1]{GZ2}:
\begin{enumerate}
\item[\bull] $A(n,q)$, generated by a grouplike element $x$, its inverse, and a $(1,x^n)$-skew primitive element $y$ satisfying $xy=qyx$, where $n\in\Znonneg$ and $q\in\kx$. (The coordinate rings $\OG$ in the third item above have the form $A(n,1)$.)
\item[\bull] $C(n)$, generated by a grouplike element $y$, its inverse, and a $(y^{n-1},1)$-skew primitive element $x$ satisfying $xy-yx= y^n-y$, where $n\in\ZZ_{\ge2}$.
\item[\bull] $B(n,p_0,\dots,p_s,q)$, generated by a grouplike element $x$, its inverse, and commuting $(1,x^{m_in})$-skew primitive elements $y_1,\dots,y_s$ satisfying $y_i^{p_i}= y_j^{p_j}$ and $xy_i= q^{m_i}y_ix$, where $n$ and $m_i= p_1p_2\cdots p_s/p_i$ are positive integers satisfying some numerical conditions we do not list here, and $q\in\kx$ is a primitive $(np_1p_2\cdots p_s/p_0)$\,th root of unity.
\end{enumerate}

All the above examples, except for the last, are Auslander-regular and Cohen-Macaulay of global dimension $2$, while the last is Auslander-Gorenstein and Cohen-Macaulay of injective dimension $2$ \cite[Proposition 0.2 and proof]{GZ2}.

A classification of (H2) was obtained under an additional homological assumption, equivalent to the existence of an infinite dimensional commutative quotient algebra \cite[Theorem 0.1, Proposition 3.8]{GZ2}.

\begin{theorem}
{\rm [Goodearl-Zhang]} If $H$ satisfies {\rm(H2)} and $\Ext^1_H({}_Hk,{}_Hk)\ne 0$, then $H$ is isomorphic to one of the examples described above.
\end{theorem}

While the homological assumption $\Ext^1_H({}_Hk,{}_Hk)\ne 0$ appears natural, it does not always hold, as shown by Wang, Zhang, and Zhuang \cite[Section 2]{WZZ1}.

\begin{construction}  \label{WZZeg}
{\rm [Wang-Zhang-Zhuang]} There exists a family of Hopf algebras satisfying {\rm(H2)} but $\Ext^1_H({}_Hk,{}_Hk)=0$.
\end{construction}

Construction \ref{WZZeg} is a modification of that for the $B(n,p_0,\dots,p_s,q)$ of \cite{GZ2}. Both families satisfy polynomial identities, and in fact, one has the following result \cite[Corollary 1.12]{WZZ1}, \cite[Theorem 0.1]{WZZ2}.

\begin{theorem}
{\rm [Wang-Zhang-Zhuang]} Any Hopf algebra $H$ satisfying {\rm(H2)} and $\Ext^1_H({}_Hk,{}_Hk)=0$ is a PI algebra.
\end{theorem}

\begin{conjecture}
{\rm [Wang-Zhang-Zhuang]} Any Hopf algebra $H$ satisfying {\rm(H2)} and $\Ext^1_H({}_Hk,{}_Hk)=0$ is isomorphic to one of those in Construction {\rm\ref{WZZeg}}.
\end{conjecture}

We end with a more ambitious problem, to line up the mentioned classifications in Gelfand-Kirillov dimensions $1$ and $2$:

\begin{problem}
Classify the prime Hopf algebras, over an algebraically closed field of characteristic zero, which are affine or noetherian with Gelfand-Kirillov dimension at most $2$.
\end{problem}

\section*{Acknowledgements} 
We thank N. Andruskiewitsch, K.A. Brown, S. Montgomery, C. Ohn, H.-J. Schneider, and J.J. Zhang for correspondence, discussions, examples, and references.


\end{document}